\documentclass[reqno]{amsart}
\usepackage{amssymb,amsmath,amsthm, setspace}
\numberwithin{equation}{section}

\newcommand{\fpdd}{\textsf{fpd}}
\newcommand{\fpd}{\textsf{fpd} }

\newcommand{\del}{\delta}

\newcommand{\R}{\mathbb{R}}
\newcommand{\C}{\mathbb{C}}
\newcommand{\N}{\mathbb{N}}
\newcommand{\Z}{\mathbb{Z}}

\renewcommand{\rho}{\varrho}
\newcommand{\hhh}{h}
\newcommand{\dpsi}{\Delta_F}
\newtheorem{thm}{Theorem}
\newtheorem{lem}{Lemma}

\newtheorem{cor}{Corollary}

\renewcommand{\mod}[1]{\hspace{-2.9mm}\pmod{#1}}
\newcommand{\x}{{\bf x}}
\newcommand{\e}{\emph}

\newcommand{\bfP}{\mathbb{P}}

\newcommand{\ma}{\mathbf}
\newcommand{\ben}{\begin{enumerate}}
\newcommand{\een}{\end{enumerate}}
\newcommand{\beq}{\begin{equation}}
\newcommand{\eeq}{\end{equation}}

\newcommand{\ve}{\varepsilon}

\newcommand{\mcal}{\mathcal}
\newcommand{\lab}{\label}
\newcommand{\al}{\alpha}

\newcommand{\D}{\Delta}
\newcommand{\be}{\beta}

\newcommand{\colt}[2]{\genfrac{}{}{0pt}{1}{#1}{#2}}

\newcommand{\eqm}[2]{\equiv #1 \pmod{#2}}

\renewcommand{\leq}{\leqslant}
\renewcommand{\geq}{\geqslant}

\DeclareMathOperator{\disc}{disc}

\theoremstyle{definition}
\newtheorem*{ack}{Acknowledgement}

\begin{document}

\title{Sums of arithmetic functions over values of binary forms}

\author{R. de la Bret\`eche}
\author{T.D. Browning}
\address{Universit\'e Paris-Sud, B\^atiment 425, 91405 Orsay Cedex}
\email{regis.de-la-breteche@math.u-psud.fr}

\address{School of Mathematics,  
University of Bristol, Bristol BS8 1TW}
\email{t.d.browning@bristol.ac.uk}


\begin{abstract}
Given a suitable arithmetic function $h:\N\rightarrow
\R_{\geq 0}$, and a binary form $F\in\Z[x_1,x_2]$, we investigate the average order of $h$ as it ranges
over the values taken by $F$. A general upper bound is obtained for
this quantity, in which the dependence upon the coefficients of $F$ is made completely explicit.
\end{abstract}
\maketitle

\section{Introduction}
 
This paper is concerned with the average order of certain
arithmetic functions,  as they range over the values taken by binary
forms. We shall say that a non-negative sub-multiplicative function
$\hhh$ belongs to the class $\mcal{M}(A,B)$ if there exists a constant $A$ such that
$\hhh(p^\ell)\leq A^\ell$ for all primes $p$ and all $\ell\in \N$, and 
there is a function $B=B(\ve)$ such that for any $\ve>0$ one has $\hhh(n)\leq B n^\ve$ for all $n\in\N$.
Let $F\in\Z[x_1,x_2]$ be a non-zero binary  form of degree $d$, such that the discriminant $\disc(F)$ is non-zero.
Such a form takes the shape
\begin{equation}\lab{19-case2}
F(x_1,x_2)=x_1^{d_1}x_2^{d_2}G(x_1,x_2),
\end{equation}
for integers $d_1,d_2\in \{0,1\}$, and a non-zero binary form
$G\in\Z[x_1,x_2]$ of degree $d-d_1-d_2$.  Moreover, we may assume that
$\disc (G)\neq 0$ and $G(1,0)G(0,1)\neq 0$.

Given a function $\hhh \in \mcal{M}(A,B)$ and a binary form $F$ as
above, the primary goal of this paper is to bound the size of the sum 
$$
S(X_1,X_2;h,F):=\sum_{1\leq n_1 \leq X_1}\sum_{1\leq n_2 \leq X_2}  \hhh (|F(n_1,n_2)|),
$$
for given $X_1,X_2>0$.  For certain choices of $h$ and $F$ it is
possible to prove an asymptotic 
formula for this quantity.  When $h=\tau$ is the
usual divisor function, for example, Greaves \cite{greaves} has shown
that there is a constant $c_F>0$ such that 
$$
S(X,X;\tau,F)=c_FX^2\log X\big(1+o(1)\big),
$$
as $X\rightarrow\infty$, when $F$ is irreducible of degree $d=3$. 
This asymptotic formula has been extended to 
irreducible quartic forms by Daniel \cite{daniel}.
When $d\geq 5$ there are no binary forms $F$ for which an
asymptotic formula is known for $S(X,X;\tau,F)$. 
In order to illustrate the main results in this article, however,
we shall derive an upper bound for $S(X,X;\tau,F)$ of the expected order of magnitude.
The primary aim of this work is to provide general upper bounds for
the sum $S(X_1,X_2;h,F)$, in which the dependence upon the coefficients of the form $F$ is made
completely explicit. We will henceforth allow the implied constant in
any estimate to depend upon the degree of the polynomial that is under
consideration. Any further dependences will be indicated by an
appropriate subscript.

Before introducing our main result, we first need to introduce some more
notation. We shall write $\|F\|$ for the maximum modulus of the
coefficients of a binary integral form $F$, and we shall say that $F$
is primitive if the greatest common divisor of its coefficients is $1$.
These definitions extend in an obvious way to all polynomials with
integer coefficients.
Given any $m\in\N$, we set
\begin{equation}\lab{rho*}
\rho^*_F(m):=\frac{1}{ \varphi(m)}\#\Big\{ (n_1,n_2)\in (0,m]^2\,:\,
\begin{array}{ll} \gcd(n_1,n_2,m)=1 \\
F(n_1,n_2)\equiv 0
\mod{m}\end{array}\Big\},
\end{equation}
where $\varphi$ is the usual Euler totient function. 
The arithmetic function $\rho^*_F$ is multiplicative, and has already played an
important role in the work of Daniel~\cite{daniel}.
Finally, we define
\begin{equation}\label{defpsi}
\psi(n):=\prod_{p\mid n}\Big(1+\frac{1}{p}\Big),
\end{equation} 
and 
\begin{equation}\label{D-F}
\dpsi:=\psi\big(\disc(F)\big).
\end{equation}
We are now ready to reveal our main result.

\begin{thm}\lab{n''}
Let $\hhh \in \mcal{M}(A,B)$, 
$\delta\in (0,1)$ and let $X_1,X_2>0$.  Let $F\in\Z[x_1,x_2]$ be a
non-zero primitive binary form of the shape \eqref{19-case2}. 
Then there exist positive constants $c=c(A,B)$
and $C=C(A,B,d,\delta)$ such that 
$$
S(X_1,X_2;h,F) \ll_{A,B,\delta}
\dpsi^c  X_1X_2 E  
$$
for $\min\{X_1,X_2\}\geq C\max\{X_1,X_2\}^{\delta d}\|F\|^\delta$,
where $\dpsi$ is given by \eqref{D-F} and  
\begin{equation}\label{defE}
\begin{split}
E:=&
\prod_{d<p\leq
\min\{X_1,X_2\}}\Big(1+\frac{\rho_G^*(p)(\hhh(p)-1)}{p}\Big)
\prod_{i=1,2}
\prod_{p\leq X_i}\Big(1+\frac{d_i(\hhh   (p)-1)}{p}\Big).
\end{split}
\end{equation} 
\end{thm}

We shall see shortly that the condition $p>d$ ensures that
$\rho_G^*(p)<p$ in \eqref{defE}.
Our initial motivation for establishing a result of the type in
Theorem \ref{n''} arose in a rather different context. It turns
out that Theorem \ref{n''} plays an important role in the authors' recent proof
of the Manin conjecture for the growth rate
of rational points of bounded height on a certain Iskovskih surface
\cite{dp4-2a1}.  The precise result that we make use of is the following,
which will be established in the subsequent section.

\begin{cor}\lab{Corn''}
Let $\hhh \in \mcal{M}(A,B)$ and let $X_1,X_2>0$.  Let $F\in\Z[x_1,x_2]$ be a
non-zero binary form of the shape \eqref{19-case2}. 
Then we have
$$
\sum_{|n_1| \leq X_1}\sum_{|n_2| \leq X_2}  \hhh (|F(n_1,n_2)|) \ll_{A,B,\ve}
\|F\|^\ve\Big(
   X_1X_2E +\max\{X_1,X_2\}^{1+\ve}\Big),
$$
for any $\ve>0$, where $E$ is given by \eqref{defE}.  
\end{cor}

An inspection of the proof of Corollary \ref{Corn''} reveals that it is possible 
to replace the term $X^{1+\ve} $ by $X(\log X)^{A^d-1}$,
where $X=\max\{X_1,X_2\}$.  Moreover, it would not be difficult to
extend the estimates in Theorem \ref{n''} and Corollary \ref{Corn''}
to the more general class of arithmetic functions 
$\mcal{M}_1(A,B,\varepsilon)$ considered by Nair and Tenenbaum
\cite{NT}.

It is now relatively straightforward to use Theorem \ref{n''} to deduce
good upper bounds for $S(X,X;h,F)$ for various well-known
multiplicative functions $h$.  For example, on taking $h=\tau$ in
Theorem \ref{n''}, and appealing to work of Daniel \cite[\S 7]{daniel}
on the behaviour of the Dirichlet series $\sum_{n=1}^\infty \rho_F^*(n)n^{-s}$, 
it is possible to deduce the following result, 
which is new for $d\geq 5$.

\begin{cor}\lab{Corn'''}
Let $F\in\Z[x_1,x_2]$ be an irreducible binary form of degree $d$. 
Then we have
$S(X,X;\tau,F) \ll_F X^2\log X.$
\end{cor}

The primary ingredient in our work is a result due to Nair
\cite{nair}.  Given an arithmetic function $\hhh \in \mcal{M}(A,B)$,  and a
suitable polynomial $f\in \Z[x]$, Nair investigates the size of the sum 
$$
T(X;\hhh,f):=\sum_{1\leq n \leq X}  \hhh (|f(n)|),
$$
for given $X>0$.  Nair's work has since been
generalised in several directions by Nair and Tenenbaum \cite{NT}.
In order to present the version of Nair's result that we shall employ,
we first need to introduce some more notation. 
Given any polynomial $f \in \Z[x]$ and any $m\in\N$,  let
$$
\rho_f(m):=\#\{n\mod{m}: ~f(n)\equiv 0 \mod{m}\}.
$$
It is well-known that $\rho_f$ is a multiplicative function.  
On recalling the definition \eqref{rho*} of $\rho_{G}^*(p)$, 
for any binary form $G\in\Z[x_1,x_2]$ and any prime $p$,
we may therefore record the equalities
\begin{equation}\lab{rh}
\rho^*_G(p)=
\left\{
\begin{array}{ll}\rho_{G(x,1)}(p) &
\mbox{if $p\nmid G(1,0)$},\\
\rho_{G(x,1)}(p)+1&
\mbox{if $p\mid G(1,0)$}.
\end{array}
\right.
\end{equation}
One may
clearly swap the roles of the first and second variables in this expression.
It follows from these equalities that $\rho_G^*(p)<p$ for any prime
$p>\deg G$, as claimed above.

Given a positive integer $d$ and a prime number $p$, we shall
denote by $\mcal{F}_p(d)$ the class of polynomials $f\in \Z[x]$
of degree $d$, which have no repeated roots and do not have $p$ as a fixed prime
divisor. Note that a polynomial has no repeated roots if
and only if its discriminant is non-zero.  Moreover, 
recall that a polynomial $f \in \Z[x]$ is said to have fixed
prime divisor $p$ if $p \mid f(n)$ for all $n \in \Z$.  
It will be convenient to abbreviate ``fixed prime divisor"
to ``{\fpdd}'' throughout this paper.
When $f$ has  degree $d$ and is primitive, then any \fpd $p$ of $f$ satisfies $p \leq
d$. Indeed, there are at most $d$ roots of $f$ modulo $p$.
We shall write 
$$
\mcal{F}(d):=\bigcap_{p}\mcal{F}_p(d).
$$
We are now ready to reveal the version of Nair's result that we shall employ.

\begin{thm}\lab{n}
Let $\hhh\in \mcal{M}(A,B)$, let $f\in \mcal{F}(d)$ and let
$\delta\in (0,1)$. Then there exists a constant $C=C(A,B,d,\delta)$ 
such that
$$
T(X;\hhh,f) \ll_{A,B,\delta} X \prod_{p \leq
  X}\Big(1-\frac{\rho_f(p)}{p}\Big) \sum_{1\leq m \leq X} \frac{\rho_f(m)\hhh(m)}{m},
$$
for $X\geq C \|f\|^\delta$.
\end{thm}

A few remarks are in order here. First and foremost this is not quite
the main result in \cite[\S 4]{nair}.  In its present form, Theorem \ref{n}
essentially amounts to a special case of a very general result due to Nair and
Tenenbaum \cite[Eqn. (2)]{NT}. 
Following our convention introduced above, the implied constant in
this estimate is completely independent of the coefficients of $f$,
depending only upon the choices of $A,B,\delta$ and $d$.  
This uniformity will prove crucial in our deduction of Theorem \ref{n''}.
Theorem \ref{n} is in fact already implicit in the
original work of Nair \cite{nair}, and is a major step on the way
towards his upper bound
\begin{equation}\lab{18-n1'}
T(X;\hhh,f) \ll_{A,B,\delta, \disc(f)} X \prod_{p \leq
  X}\Big(1-\frac{\rho_f(p)}{p}\Big) \exp\Big(\sum_{p \leq X} \frac{\hhh(p)\rho_f(p)}{p}\Big),
\end{equation}
for $X\geq C \|f\|^\delta$.  As indicated, there is now an implicit
dependence upon the discriminant of the polynomial $f$. This arises
in passing from the term $\sum_{m} \frac{\hhh(m)\rho_f(m)}{m}$ to the term 
$\exp\big(\sum_{p} \frac{\hhh(p)\rho_f(p)}{p}\big)$.

We take this opportunity to correct an apparent oversight in recent
work of Heath-Brown \cite{h-b03}. Here, a special case of Nair's result is
used \cite[Lemma~4.1]{h-b03}, in which the dependence of the
implied constant upon the polynomial's discriminant does not seem to
have been accurately recorded.  This leads to problems in the proof of
\cite[Lemma 4.2]{h-b03}, and in particular the estimation of the sum
$S_0(m)$, since the relevant discriminant will now vary with the choice $m$.
Similar remarks apply to the estimation of $S(d,d')$ in \cite[Lemma 6.1]{h-b03}.
The proof of these two estimates can now be easily repaired: the first
by appealing to Theorem \ref{n} instead of \eqref{18-n1'}, and the
second via a straightforward application of Theorem~\ref{n''}.

\begin{ack}
Part of this work was undertaken while the second author was visiting 
the Universit\'e de Paris-Sud, the hospitality and financial support
of which is gratefully acknowledged.
\end{ack}

\section{Preliminaries}\lab{pre}

We begin this section by establishing Corollary \ref{Corn''}.
Now it is trivial to see that $\D_F\ll_\ve \|F\|^\ve$, since the
discriminant of a form can always be bounded in terms of the maximum
modulus of its coefficients, and we have 
$$
\psi(n)\leq 2^{\omega(n)}\ll_\ve n^\ve,
$$
for any $\ve>0$.  Moreover, it will suffice to establish the result
under the assumption that $F$ is primitive. Indeed, if 
$k$ is a common factor of the coefficients of $F$, then it may
extracted and absorbed into the factor $\|F\|^\ve$, since 
$h(ab)\ll_{B,\ve} a^\ve h(b)$ for $h\in \mcal{M}(A,B)$.
Let us take $\delta=\ve$ in the statement of Theorem \ref{n''}.
Suppose first that 
$\min\{X_1,X_2\}\leq C\max\{X_1,X_2\}^{d \ve}\|F\|^\ve$.
Then
since $E\ll_\ve (X_1X_2)^\ve$ in \eqref{defE}, we easily deduce that 
$$
S(X_1,X_2;h,F) \ll_{A,B,\ve}
\|F\|^\ve  \max\{X_1,X_2\}^{1+\ve}.
$$
This is satisfactory for Corollary \ref{Corn''}. In the alternative
case, Theorem \ref{n''} gives a satisfactory contribution from those
$\ma{n}$ for which $n_1n_2\neq 0$. The contribution from $n_1=0$ is 
$$
\leq 
\sum_{|n_2|\leq X_2} \hhh   (|F(0,n_2)|)
\leq  \hhh   (|F(0,1)|) 
\sum_{|n_2|\leq X_2} \hhh(n_2^{d})  
\ll_{B,\ve} \|F\|^\ve X_2^{1+\ve},
$$
since $\hhh\in \mcal{M}(A,B)$, which is also satisfactory. On arguing
similarly for the contribution from $n_2=0$, we therefore complete the
proof of Corollary \ref{Corn''}.

We now collect together the preliminary facts that we shall need in our
proof of Theorem \ref{n''}.  
Let $F\in\Z[\x]$ be a non-zero binary form of
degree $d$.  Here, as throughout our work, any boldface lowercase letter $\x$ will
mean an ordered pair $(x_1,x_2)$.
If $[\alpha_1,\be_1], \ldots, [\alpha_d,\be_d]\in
\bfP^1(\C)$ are the $d$ roots of $F$ in $\C$, then the discriminant of
$F$ is defined to be 
$$
\disc(F):=\prod_{1\leq i<j \leq d}(\al_i\be_j - \al_j\be_i)^2.
$$
It will be convenient to record the following well-known result.

\begin{lem}\lab{disc}
Let $\ma{M}\in GL_2(\Z)$. Then we have
$$
\disc(F(\ma{M}\x))=\det(\ma{M})^{d(d-1)}\disc(F).
$$
\end{lem}

We shall also require good upper
bounds for the quantity $\rho_f(p^\ell)$, for any primitive polynomial
$f\in \Z[x]$ and any prime power $p^\ell$. 
The following result may be found in unpublished work of Stefan
Daniel, the proof of which we provide here for the sake of completeness.

\begin{lem}\lab{dan}
Let $d\in \N$, let $p$ be a prime, and let $f\in \Z[x]$ be a
polynomial of degree $d$ such that $p$ does not divide all of the
coefficients of $f$.  Then we have
$$
\rho_f(p^{\ell}) \leq \min\big\{dp^{\ell-1}, 2d^3 p^{(1-1/d)\ell}\big\},
$$
for any $\ell\in\N$.
\end{lem}

\begin{proof}
The upper bound $\rho_f(p^{\ell}) \leq dp^{\ell-1}$  is trivial.
The second inequality is easy when $d=1$, or when $p$ divides all of
the coefficients of $f$ apart from the constant term, in which case
$\rho_f(p^\ell)=0$.  Thus we may proceed under the assumption that
$d\geq 2$ and $p$ does not divide all of the coefficients
in the non-constant terms.  We have
\begin{align*}
\rho_f(p^\ell)
&=\frac{1}{p^{\ell}}\sum_{a \mod{p^\ell}} \sum_{b \mod{p^\ell}}
e_{p^{\ell}}(af(b))
=\sum_{j=0}^\ell  \frac{1}{p^{j}}
\sum_{\colt{a \mod{p^j}}{p\nmid a}} \sum_{b \mod{p^j}} e_{p^{j}}(af(b)),
\end{align*}
where $e_q(z)=e^{2\pi i z/q}$, as usual.  But then 
the proof of \cite[Theorem~7.1]{vaughan} implies that each inner sum is bounded by $d^3
p^{(1-1/d)j}$ in modulus, when $j\geq 1$.  Hence 
$$
\rho_f(p^\ell)
\leq 1+ d^3(1-p^{-1})\sum_{j=1}^\ell
p^{(1-1/d)j} \leq d^3p^{(1-1/d)\ell} \frac{1-p^{-1}}{1-p^{1/d-1}}.
$$
The result then follows, since $d\geq 2$ by assumption.
\end{proof}

The remainder of this section concerns the class of primitive
polynomials $f\in\Z[x]$ which have a \fpdd. The following result is self-evident.

\begin{lem}\lab{Reduc1'} 
Let $p$ be a prime number and let $f\in\Z[x]$ be a primitive
polynomial which has $p$ as a \fpdd.
Then there exists an integer $e\geq 0$ and polynomials $q,r\in\Z[x]$, such that
\begin{equation}\lab{PQR}
f(x)=(x^p-x)q(x)+pr(x),
\end{equation}
where $q(x)=\sum_{j=0}^e a_j x^j$ for integers $0\leq a_j<p$
such that $a_{e}\neq 0$.
\end{lem}

Our next result examines the effect of making the change of
variables $x\mapsto px+k$, for integers $0\leq k< p$.

\begin{lem}\lab{Reduc1} 
Let $p$ be a prime number and let $f\in\Z[x]$ be a primitive
polynomial of the shape \eqref{PQR}.
Then for each $0\leq k<p$, there exists $\nu_k\in \Z$ such that:
\begin{enumerate}
\item
$0\leq \nu_k\leq e$.
\item
$f_{k}(x)=p^{-\nu_k-1}{f(px+k)}\in\Z[x]$ is a primitive polynomial.
\item
Suppose that $f_k$ has $p$ as a \fpdd, and is written in the form
\eqref{PQR} for suitable polynomials $q_k,r_k$. Then $e\geq p-1$ and $\deg (q_k)\leq e-p+1$.
\end{enumerate}
\end{lem}
\goodbreak

\begin{proof}
Without loss of generality we may assume that $k=0$. 
Consider the identity
$$
\frac{f(px )}{p }=x(p^{p-1}x^{p-1}-1)q(px)+r(px),
$$
and let $b_j$ be the  $j$-th coefficient of $r(x)$.
It is not hard to see that the coefficient of $x^{j+1}$ in $f(px)/p$
is equal to
\begin{equation}\lab{15-bi}
(a_{j-p+1}-a_j)p^{j}+b_{j+1}p^{j+1},
\end{equation}
where we have introduced the convention that $a_j=0$ for
each negative index $j$.
Let $\nu_0$
be the $p$-adic order of the greatest common divisor of the
coefficients of the polynomial $f(px)/p$, and write
$$
f_0(x)=\frac{f(px)}{p^{\nu_0+1}}.
$$
It is clear that $f_0$ is a primitive polynomial with integer coefficients.
Moreover, if $e_0$ denotes the smallest index $j$ for which $a_j\neq
0$ in $q(x)$, then it is not hard to deduce from \eqref{15-bi} that 
$\nu_0\leq e_0$.
In particular we have $0\leq \nu_0 \leq e$. This is enough to
establish the first two parts of the lemma. 

It remains to consider the possibility that $f_0$ has $p$ as a \fpdd. 
Suppose first that $\nu_0< e_0$.  Then $f_0(x)\eqm{g_0(x)}{p}$,
with 
$$
g_0(x)=\sum_{\ell=0}^{\nu_0} b_\ell p^{\ell-\nu_0}x^\ell.
$$
If $g_0$ has $p$ as a \fpdd, then one may write it in
the form \eqref{PQR} for suitable $q_0,r_0 \in\Z[x]$. But then 
$$
0\leq \deg(q_0)\leq \nu_0-p< e_0-p\leq e-p,
$$
which is satisfactory for the final part of the lemma. Suppose now that
$\nu_0=e_0.$  Then $f_0(x)\eqm{g_0(x)}{p}$, with 
$$
g_0(x)= -a_{e_0}x^{e_0+1}+\sum_{\ell=0}^{e_0}
b_\ell p^{\ell-\nu_0} x^\ell.
$$ 
Arguing as above, if $g_0$ has $p$ as a \fpdd, then one may 
write it in the form \eqref{PQR} for suitable $q_0,r_0\in\Z[x]$ such that
$$
0\leq \deg(q_0)= e_0+1-p\leq e+1-p.
$$
This therefore completes the proof of Lemma \ref{Reduc1}.
\end{proof}

Our final result combines Lemmas \ref{Reduc1'} and \ref{Reduc1} in
order to show that there is always a linear
change of variables that takes a polynomial with \fpd
$p$ into a polynomial which doesn't have $p$ as a \fpdd.

\begin{lem}\lab{Reduc2}
Suppose that $f\in\Z[x]$ is a primitive polynomial that takes the
shape \eqref{PQR} and has non-zero discriminant. Then there exists a
non-negative integer $\del\leq e$,
and positive integers $\mu_0,\ldots,\mu_\del$ with
\begin{equation}\lab{14-tdog}
\mu_0+\cdots+\mu_\del \leq (e+1)^2,
\end{equation} 
such that the polynomial
\begin{equation}\lab{14-scat}
g_{k_0,\ldots,k_\del}(x)=\frac{f(p^{\del+1}x+p^\del k_\del+\cdots+pk_1+k_0)}{p^{\mu_0+\cdots+\mu_\del}}
\end{equation}
belongs to $\mcal{F}_p(d)$,
for any $k_0,\ldots,k_\del\in\Z\cap[0,p)$.
\end{lem}

\begin{proof}
Our argument will be by induction on the degree $e$ of $q$.
We begin by noting that the degree of $f$ is preserved under any
linear transformation of the shape $x\mapsto ax+b$, provided that
$a\neq 0$.  Similarly, in view of Lemma \ref{disc}, 
the discriminant will not vanish under any such transformation.
Thus it suffices to show that there exists a non-negative integer
$\del\leq e$, and positive integers $\mu_0,\ldots,\mu_\del$, such that
\eqref{14-tdog} holds and the polynomial \eqref{14-scat} has
integer coefficients but doesn't have $p$ as a \fpdd.

Let $k_0$ be any integer in the range
$0\leq k_0<p$.  Then it follows from Lemma~\ref{Reduc1} that 
there exists $\nu_0\in \Z$ such that $0\leq \nu_0\leq e$ and 
$$
f_{k_0}(x)=p^{-\nu_0-1}{f(px+k_0)}
$$
is a primitive polynomial with integer coefficients.
If $e<p-1$ then the final part of this result implies that 
$f_{k_0}$ does not contain $p$ as a \fpdd, and so must belong to $\mcal{F}_p(d)$.
In this case, therefore, the statement of Lemma \ref{Reduc2} holds with $\del=0$,
$\mu_0=\nu_0+1$ and $g_{k_0}=f_{k_0}$. 
This clearly takes care of the inductive base $e=0$, since then
$\del=0$ and $\mu_0=1$.
Suppose now that $e\geq p-1$ and $f_{k_0}$ contains $p$ as a \fpdd.
Then $f_{k_0}$ can be written in the form
\eqref{PQR} for suitable polynomials $q',r'$ such that $\deg(q')=e'\leq e-p+1$.
We may therefore apply the inductive hypothesis to
conclude that there exists a non-negative integer $\del'\leq e'$,
and positive integers $\mu_0',\ldots,\mu_{\del'}'$ with
\begin{equation}\lab{14-tdog'}
\mu_0'+\cdots+\mu_{\del'}' \leq  (e'+1)^2,
\end{equation} 
such that the polynomial
$$
\frac{f_{k_0}(p^{\del'+1}x+p^{\del'}
  k_{\del'}'+\cdots+p{k_1}'+{k_0'})}{p^{\mu_0'+\cdots+\mu_{\del'}'}}=
\frac{f(p^{\del'+2}x+p^{\del'+1} k_{\del'}'+\cdots+p{k_0'}+k_0)}{p^{\mu_0'+\cdots+\mu_{\del'}'+\nu_0+1}}
$$
belongs to $\mcal{F}_p(d)$,
for any $k_0,k_0',\ldots,k_{\del'}'\in\Z\cap[0,p)$.
Let $\del=\del'+1$, let $k_i'=k_{i+1}$ for $i\geq 0$,
and write
$$
\mu_0=\nu_0+1, \quad \mu_i=\mu_{i-1}', 
$$
for $i\geq 1$.  Then it
follows that 
$g_{k_0,\ldots,k_\del}(x)\in \mcal{F}_p(d)$, in the
notation of \eqref{14-scat}, for any $k_0,\ldots,k_{\del}\in\Z\cap[0,p)$.
Moreover, we clearly have $\del\leq e'+1\leq e-p+2 \leq e,$ 
and \eqref{14-tdog'} gives
$$
\mu_0+\cdots+\mu_\del 
\leq  (e-p+2)^2+(e+1)\leq e^2+e+1\leq (e+1)^2.
$$
Thus \eqref{14-tdog} also holds, which therefore 
completes the proof of Lemma \ref{Reduc2}.
\end{proof}

Suppose that $f\in \Z[x]$ is a primitive polynomial that
takes the shape \eqref{PQR} for some prime $p$, but 
which does not have $q$ as a \fpd for any prime $q<p$.  
Then for any $a\in \Z$, the linear polynomial $p^{\del+1}x+a$ runs
over a complete set of residue classes modulo $q$ as $x$ does. Thus it follows
from the statement of Lemma \ref{Reduc2} that
$$
g_{k_0,\ldots,k_\del}\in \bigcap_{q\leq p} \mcal{F}_q(d),
$$
for any $k_0,\ldots,k_\del\in\Z\cap[0,p)$, where the intersection is over all primes $q\leq p$.

\section{Proof of Theorem \ref{n''}} 

We are now ready to proceed with the proof  of Theorem
\ref{n''}. Suppose that $ X_2\geq X_1\geq 1$,
and let $F\in\Z[\x]$ be a primitive form of the shape \eqref{19-case2},
Let $d'=d-d_2$ and $d''=d-d_1-d_2$. We may therefore write
$$
G(\x)=\sum_{j=0}^{d''}a_j x_1^{d''-j}x_2^{j},
$$
for $a_j\in\Z$ such that $\gcd(a_0,\ldots,a_{d''})=1$ and $a_0a_{d''}\neq 0$.
We begin this section by recording the following easy result.
 
\begin{lem}\lab{disc2}
Let $p$ be a prime. Then we have $p \mid \disc(F)$ for any $p\mid \gcd (a_0,a_1).$
Moreover, if $d_2=1$, then 
we have $p \mid \disc(F)$ for any $p\mid a_0.$
\end{lem}
\begin{proof} 
The first fact follows on observing that the reduction of $F$ modulo
$p$ has $x_2^2$ as a factor if $p\mid \gcd(a_0,a_1)$. If $d_2=1$, then
the same conclusion holds provided only that $p\mid a_0$.
The statement of the lemma is now obvious.
\end{proof}

We intend to apply Theorem \ref{n}, for which we shall fix one of
the variables at the outset.
Let $q_m:=\gcd (a_0 ,a_1m,\ldots, a_{d''}m^{d''})$, for any $m\in\N$,  and define 
$$
f_{n_2}(x):=\frac{x^{d_1}G(x,n_2)}{q_{n_2}}.
$$  
Then it is clear that $f_{n_2}$ is a primitive polynomial of degree
$d'$ with integer coefficients.  Moreover, we have
\begin{equation}\lab{19-step1}
S(X_1,X_2;h,F) \leq
\sum_{1\leq n_2  \leq X_2} \hhh\big(n_2^{d_2}q_{n_2}\big) 
\Big|\sum_{1\leq n_1 \leq X_1}\hhh\big(|f_{{n_2}}(n_1)|\big)\Big|.
\end{equation}
We now want to apply Theorem \ref{n} to estimate the inner sum.
For this we must deal with the possibility that $f_{m}$ contains a \fpdd.
Since $f_m$ is primitive of degree $d'$,
the only possible {\fpdd}s are the primes $p\leq d'$.

Suppose that $f_m$ has $p_1<\cdots<p_r$ as {\fpdd}s. We shall
combine a repeated application of Lemma \ref{Reduc2} with the
observation made at the close of \S \ref{pre}. 
This leads us to the conclusion that 
there exist non-negative integers $\del_1,\ldots,\del_r\leq d-2$, together
with positive integers $m_1,\ldots,m_r\leq d^2$, such that
$$
g_\beta(x):=\frac{f_{n_2}(p_1^{\del_1+1}\cdots
  p_r^{\del_r+1}x+\beta)}{p_1^{m_1}\cdots p_r^{m_r}} \in \mcal{F}(d'),
$$
for any $\beta$ modulo $p_1^{\del_1+1}\cdots p_r^{\del_r+1}$.
It will be convenient to write 
$$
\alpha:=p_1^{\del_1+1}\cdots p_r^{\del_r+1}, \quad \gamma:=p_1^{m_1}\cdots p_r^{m_r}.
$$
Then it follows from Lemma \ref{disc} that
\begin{equation}\lab{20-hdisc}
\begin{split}
\disc(g_\beta)=
\disc\Big( \frac{(\al x+\beta)^{d_1}G(\alpha x+\beta,n_2)}{\gamma q_{n_2}}\Big)
&=
\disc\Big( \frac{F(\alpha x+\beta,n_2)}{\gamma
  q_{n_2}n_2^{d_2}}\Big)\\
&=\Big(\frac{\alpha^{d}n_2^{d-2d_2}}{\gamma^{2}q_{n_2}^2}\Big)^{d-1}\disc(F).
\end{split}
\end{equation}
Note that
$\alpha \leq  d^{r(d-1)} \leq d^{d^2}$ and $\gamma\leq d^{rd^2}\leq
d^{d^3}.$
In particular there are just $O(1)$ choices for $\beta$ modulo $\alpha$, and
$\hhh(\gamma)\ll_{B} 1$.

Our investigation so far has therefore led us to the inequality
\begin{equation}\lab{20-conc}
\sum_{1\leq n_1 \leq X_1} \hhh(|f_{n_2}(n_1)|) \ll_{B}
\sum_{\alpha} \sum_{\beta \mod{\alpha}}
\sum_{1\leq n_1 \leq X_1 } \hhh(|g_\beta(n_1)|),
\end{equation}
in \eqref{19-step1}, with $g_\beta \in \mcal{F}(d')$.  
It will now suffice to
apply Theorem \ref{n} to estimate the inner sum, which we henceforth
denote by $U(X_1)$.  Note that $\|g_\beta\|\ll \|f_{n_2}\| \ll {n_2}^d\|F\| \leq X_2^d\|F\|$. Hence
it follows from Theorem \ref{n} that for any $\delta\in(0,1)$ we have
\begin{equation}\lab{20-app}
U(X) \ll_{A,B,\delta} X \prod_{p \leq
  X}\Big(1-\frac{\rho_{g_\beta}(p)}{p}\Big) \sum_{1\leq m \leq X}
  \frac{\rho_{g_\beta}(m)\hhh(m)}{m}, 
\end{equation}
for $X\gg_{A,B,\delta} X_2^{\delta d}\|F\|^\delta$.
In estimating the right hand side of \eqref{20-app}, we shall find
that the result is largely independent of the choice of $\beta$. To
simplify our exposition, therefore, it will be convenient to write
$g=g_\beta$ in what follows.

We begin by estimating the sum over $m$ that appears in \eqref{20-app}.
On combining the sub-multiplicativity of $\hhh$ with the
multiplicativity of $\rho_g$, we see that
\begin{equation}\lab{20-split}
 \sum_{1\leq m \leq X} \frac{\rho_{g}(m)\hhh(m)}{m} \leq
\prod_{p\leq X}\Big(1+\frac{\rho_{g}(p)\hhh(p)}{p} 
+\sum_{\ell\geq 2} \frac{\rho_{g}(p^\ell)\hhh(p^\ell)}{p^\ell}\Big).
\end{equation}
We must therefore examine the behaviour of the function $\rho_g  (p^\ell)$
at various prime powers $p^\ell$.  
This is a rather classic topic and
the facts that we shall use may all be found in the book of Nagell \cite{nagell},
for example.
Now an application of Lemma \ref{dan} reveals that
$$
\rho_g  (p^{\ell}) \leq \min\big\{d'p^{\ell-1}, 2{d'}^3 p^{(1-1/d')\ell}\big\},
$$
for any $\ell\in\N$, since $p$ does not divide all of the coefficients
of $g$.  Moreover, it is well-known that
$$
\rho_g  (p^\ell) \leq d',
$$
if $p\nmid \disc(g  )$ or if $\ell=1$.
In view of the fact that 
$\hhh(p^\ell)\leq \min\{A^\ell,Bp^{\ell\ve}\},$ 
for any $\ve>0$, we therefore deduce that
\begin{align*}
\sum_{\ell\geq 1}
\frac{\rho_g  (p^\ell)\hhh(p^\ell)}{p^\ell} 
&\leq 
d'\sum_{1\leq \ell\leq d}
\frac{\hhh(p^\ell)p^{\ell-1}}{p^\ell} 
+
2{d'}^3\sum_{\ell>d} \frac{\hhh(p^\ell)p^{(1-1/d')\ell}}{p^\ell} 
\ll_{A,B}\frac{1}{p}. 
\end{align*}
for any prime $p\mid \disc(g  )$.  When $p\nmid \disc(g)$ we obtain 
\begin{align*}
\sum_{\ell\geq 2}
\frac{\rho_g  (p^\ell)\hhh(p^\ell)}{p^\ell} 
&\leq 
d'\sum_{\ell\geq 2}
\frac{\hhh(p^\ell)}{p^\ell} 
\ll_{B,\ve} p^{-2(1-\ve)}.
\end{align*}
Now \eqref{20-hdisc} implies that 
$
\psi\big(\disc(g  )\big)\leq \psi\big(\alpha n_2 \disc(F)\big)
\ll \dpsi \psi(n_2),
$
where $\psi$ is given by \eqref{defpsi} and $\dpsi$ is given by \eqref{D-F}.
Drawing our arguments together, therefore, we have so far shown that
there is a constant $c_1=c_1(A,B)$ such that 
$$
 \sum_{1\leq m \leq X} \frac{\rho_g  (m)\hhh(m)}{m} \ll_{A,B}
\dpsi^{c_1}\psi(n_2)^{c_1}\prod_{\colt{d<p\leq X}{p\nmid
\disc(g)}}\Big(1+\frac{\rho_g  (p)\hhh(p)}{p}\Big),
$$
in \eqref{20-split}. Suppose now that $p>d>d'.$
Then one has $\rho_g  (p)=\rho_{f_{n_2}}(p)$. 
We claim that $p\nmid q_{n_2}$ provided that $p\nmid n_2\disc(F)$.
But this follows immediately from the fact that $\gcd(a_0,\ldots,a_{d''})=1$.
Hence we have
\begin{equation}\lab{capote}
\rho_g  (p)=\rho_{x^{d_1}G(x,n_2)}(p)=\rho_{x^{d_1}G(x,1)}(p)=
\rho_{ G(x,1)}(p)+d_1,
\end{equation}
provided that $p>d$ and $p\nmid n_2\disc(F)$.
We may therefore conclude that there is a constant $c_2=c_2(A,B)>c_1$ such that 
\begin{align*}
 \sum_{1\leq m \leq X} 
\hspace{-0.1cm}
\frac{\rho_{g }(m)\hhh(m)}{m}   \ll_{A,B}&
\dpsi^{ c_2}\psi(n_2)^{ c_2}
\prod_{d<p\leq X}
\hspace{-0.2cm}
\Big(1+\frac{\rho_{G(x,1)}(p)\hhh(p)}{p}\Big)
\prod_{p\leq X}
\hspace{-0.1cm}
\Big(1+\frac{d_1\hhh(p) }{p}\Big).
\end{align*}

We now turn to the size of the product over $p$ that appears in 
\eqref{20-app}, for which we shall use the relation \eqref{capote} for
any prime $p$ such that $p>d$ and $p\nmid n_2\disc(F)$.
Thus there is a constant $c_3=c_3(A,B)$ such that
\begin{align*} 
\prod_{p \leq X} \Big(1-\frac{\rho_{g}(p)}{p}\Big)
&\ll 
\prod_{\colt{d<p \leq X}{p\nmid n_2\disc(F)}}\Big(1-\frac{\rho_{G(x,1)}(p)}{p}\Big)
\prod_{\colt{p \leq X}{p\nmid n_2\disc(F)}}\Big(1-\frac{d_1}{p}\Big) \\
&\ll  \dpsi^{c_3}\psi(n_2)^{c_3}
\prod_{d<p \leq X}\Big(1-\frac{\rho_{G(x,1)}(p)}{p}\Big)\prod_{p \leq
  X} \Big(1-\frac{d_1}{p}\Big).
\end{align*}
Let
$$
E_1:=\prod_{d<p \leq X_1}\Big(1-\frac{\rho_{G(x,1)}(p)}{p}\Big)
\Big(1+\frac{\rho_{G(x,1)}(p)\hhh(p)}{p}\Big) 
\prod_{p \leq
  X_1}\Big(1-\frac{d_1}{p}\Big)\Big(1+\frac{d_1\hhh(p)}{p}\Big),
$$
and set $c_4=c_2+c_3$. Then we have shown that
$$
U(X_1) \ll_{A,B,\delta} \dpsi^{c_4}\psi(n_2)^{c_4} X_1 E_1
$$
in \eqref{20-app}, provided that $X_1\gg_{A,B,\delta} X_2^{\delta
d}\|F\|^\delta$.  This latter inequality holds by the assumption made in the statement
of Theorem \ref{n''}.

Once substituted into 
\eqref{19-step1} and \eqref{20-conc}, we may conclude that
\begin{equation}\label{maj36}
S(X_1,X_2;h,F) \ll_{A,B,\del }\dpsi^{c_4}X_1 E_1 V_{d_2}(X_2),
\end{equation}
where
$$
V_{d_2}(X_2)=\sum_{1\leq n_2  \leq X_2}\psi(n_2)^{c_4} \hhh(n_2^{d_2}q_{n_2}). 
$$
We shall estimate $V_0(X_2)$ and $V_1(X_2)$ with a further application
of Theorem \ref{n}.  To begin with we note that for any prime $p$ we have
$$
q_p=\left\{\begin{array}{ll}
p,& \mbox{if $p\mid a_0$ and $p\nmid a_1$,}\\
1,& \mbox{if $p\nmid a_0$.} 
\end{array}
\right.
$$
When $p^2\mid a_0$ and $p\mid a_1$ it is clear that $q_p$ has $p^2$ as
a factor.  Lemma \ref{disc2} implies that this can only happen
when $p\mid \disc(F)$.

Suppose first that $d_2=0$. Then the arithmetic function $n\mapsto
\psi(n)^{c_4} \hhh( q_n)$ 
satisfies the conditions of Theorem \ref{n}. Applying this
result with the polynomial $f(x)=x$, as we clearly may, it therefore
follows that there is a constant $c_5=c_5(A)$ such that
\begin{align*} 
V_{0}(X_2)
&\ll_A {X_2}\prod_{p\leq X_2}
\Big(1- \frac{1}{p}\Big) \prod_{\colt{p\leq X_2}{p\nmid a_0}}
\Big(1+ \frac{1}{p}\Big) 
\prod_{\colt{p\leq X_2}{p\mid a_0}}\Big(1+ \frac{\hhh(q_p)}{p}\Big)\\
& \ll_A \dpsi^{c_5}X_2
\prod_{p\leq X_2}\Big(1- \frac{1}{p}\Big) 
\prod_{\colt{p\leq X_2}{p\nmid a_0}}
\Big(1+ \frac{1}{p}\Big) 
\prod_{\colt{p\leq X_2}{p\mid a_0}}\Big(1+ \frac{\hhh(p)}{p}\Big)\\
& \ll_A \dpsi^{c_5}X_2
\prod_{\colt{p\leq X_2}{p\mid a_0}}\Big(1- \frac{1}{p}\Big) 
\prod_{\colt{p\leq X_2}{p\mid a_0}}\Big(1+ \frac{\hhh(p)}{p}\Big),
\end{align*}
for $X_2\gg_{A,B} 1$.
Recall the identities \eqref{rh}.  Then on inserting this bound into \eqref{maj36}, we therefore obtain the expected
bound in Theorem \ref{n''} since
\begin{align*}
\prod_{\colt{p\leq X_2}{p\mid a_0}}\Big(1- \frac{1}{p}\Big) 
\prod_{d<p\leq X_1}\Big(1-\frac{\rho_{G(x,1)}(p)}{ p}\Big)
\ll \prod_{d<p \leq X_1}\Big(1-\frac{\rho^*_{G}(p) }{p}\Big),
\end{align*}
and 
\begin{align*}
\prod_{\colt{p\leq X_2}{p\mid a_0}}\Big(1+ \frac{\hhh(p)}{p}\Big)
\prod_{d<p \leq
X_1}\Big(1+\frac{\rho_{G(x,1)}(p)\hhh(p)}{p}\Big)
\ll_\del\prod_{d<p \leq
X_1}\Big(1+\frac{\rho^*_{G}(p)\hhh(p)}{p}\Big).
\end{align*}
Here we have used the elementary fact that there are at most
$\delta^{-1}$ primes $p$ such that $p\mid a_0$ and $p>a_0^{\del}$.

Let us now turn to the case $d_2=1$.  In particular it follows from 
Lemma \ref{disc2} that  $p\mid \disc(F)$ when $p\mid a_0$. Now the
function $n\mapsto \psi(n)^{c_4} \hhh(n q_n)$ again 
verifies the conditions of Theorem \ref{n}. Thus we deduce that there
exists a constant $c_6=c_6(A)$ such that
\begin{align*} 
V_1(X_2)
&\ll_A \dpsi^{c_6} X_2\prod_{p \leq
  X_2}\Big(1-\frac{d_2}{p}\Big)\prod_{\colt{p \leq
  X_2}{p\nmid a_0}}\Big(1+\frac{d_2\hhh(p)}{p}\Big)\\
&\leq \dpsi^{c_6} X_2\prod_{p \leq
  X_2}\Big(1-\frac{d_2}{p}\Big)\Big(1+\frac{d_2\hhh(p)}{p}\Big).
\end{align*}
On inserting this into \eqref{maj36}, we easily derive the desired upper bound.
This completes the proof of Theorem \ref{n''} when $X_2\geq
X_1\geq 1$. The treatment of the case in which
$X_1\geq X_2\geq 1$ is handled in precisely the same way, by 
changing the order of summation at the outset.


\begin{thebibliography}{99}

\bibitem{dp4-2a1}
R. de la Bret\`eche and T.D. Browning,
Manin's conjecture for a certain Iskovskih surface. {\em In
  preparation}, 2006.

\bibitem{daniel}
S. Daniel,
On the divisor-sum problem for binary forms.
{\em J. Reine Angew. Math.} {\bf 507} (1999), 107--129.

\bibitem{greaves}
G. Greaves,
On the divisor-sum problem for binary cubic forms.
{\em Acta Arith.} {\bf 17} (1970), 1--28.

\bibitem{h-b03}
D.R. Heath-Brown,
Linear relations amongst sums of two squares.
{\em Number theory and algebraic geometry}, 133--176,
London Math. Soc. Lecture Note Ser. {\bf 303}
CUP, 2003. 

\bibitem{nagell} T. Nagell,
{\em Introduction to number theory}. 2nd ed., Chelsea, 1964.

\bibitem{nair} M. Nair,
Multiplicative functions of polynomial values in short intervals.
{\em Acta Arith.} {\bf 62} (1992), 257--269.

\bibitem{NT} M. Nair and G. Tenenbaum,
Short sums of certain arithmetic functions. {\em Acta Math.} {\bf 180} (1998),
119--144.

\bibitem{vaughan}
R.C. Vaughan,
\e{The Hardy--Littlewood method}.
2nd ed., CUP, 1997.




\end{thebibliography}
\end{document}